\documentclass[a4paper]{amsart}

\usepackage{graphicx}
\usepackage[T1]{fontenc}
\usepackage{amssymb}
\usepackage{dsfont}
\numberwithin{equation}{section}
\newtheorem{thm}{Theorem}[section]
\newtheorem{lemma}[thm]{Lemma}
\newtheorem{cor}[thm]{Corollary}
\newtheorem{dfn}{Definition}[section]
\newtheorem{exm}{Example}[section]

\title{Density estimation with quadratic loss: a confidence intervals method}
\author{Pierre Alquier}
\address{Laboratoire de Probabilités et Modèles Aléatoires
\\
and Laboratoire de Statistique, CREST
\\
3, avenue Pierre Larousse
\\
92240 Malakoff, France.
}
\urladdr{http://www.crest.fr/pageperso/alquier/alquier.htm}
\date{\today}
\email{alquier@ensae.fr}
\thanks{I would like to thank my PhD advisor, Professor Olivier Catoni, for his kind and constant help,
and Professors Patrice Bertail, Emmanuelle Gautherat and Hugo Harari-Kermadec for their remark that 
Panchenko's lemma could improve theorem \ref{firstbound}.}

\begin{document}

\begin{abstract}
In \cite{iterative}, a least square regression estimation procedure was proposed: first, we condiser
a family of functions $f_{k}$ and study the properties of an estimator in every unidimensionnal model
$\{\alpha f_{k},\alpha\in\mathds{R}\}$;
we then show how to aggregate these estimators. The purpose of this paper is to extend this method to the case
of density estimation.
We first give a general overview of the method, adapted to the density estimation problem. We then show that this leads to adaptative estimators,
that means that the estimator reaches the best possible rate of convergence (up to a $\log$ factor).
Finally we show some ways to improve and generalize the method.
\end{abstract}

\maketitle

\section{Introduction: the density estimation setting}

Let us assume that we are given a measure space $(\mathcal{X},\mathcal{B},\lambda)$
where $\lambda$ is positive and $\sigma$-finite,
and a probability measure
$P$ on $(\mathcal{X},\mathcal{B})$ such that $P$ has a density with respect to $\lambda$:
$$ P(dx) = f(x) \lambda(dx). $$
We assume that we observe a realisation of the canonical process $(X_{1},...,X_{N})$
on $(\mathcal{X}^{N},\mathcal{B}^{\otimes N},P^{\otimes N})$.
Our objective here is to estimate $f$ on the basis of the observations $X_{1},...,X_{N}$.

More precisely, let $\mathcal{L}^{2}(\mathcal{X},\lambda)$ denote the set of all measurables functions from $(\mathcal{X},\mathcal{B})$
to $(\mathds{R},\mathcal{B}_{R})$ where $\mathcal{B}_{R}$ is the Borel $\sigma$-algebra on $\mathds{R}$. We will write
$ \mathcal{L}^{2}(\mathcal{X},\lambda)=\mathcal{L}^{2}$ for short.
Remark that $f\in\mathcal{L}^{2}$.
Let us put, for any $(g,h)\in\left(\mathcal{L}^{2}\right)^{2}$:
$$ d^{2}(g,h)=\int_{\mathcal{X}} \Bigl(g(x)-h(x)\Bigr)^{2}\lambda(dx),$$
and let $\|.\|$ and $\left<.,.\right>$ denote the corresponding norm and scalar product.
We are here looking for an estimator $\hat{f}$ that tries to minimize our objective:
$$ d^{2}(\hat{f},f). $$

Let us choose an integer $m\in\mathds{N}$ and a family of functions $(f_{1},...,f_{m})\in\left(\mathcal{L}^{2}\right)^{m}$.
There is no particular asumptions about this family: it is not necessarily linearly independant for example.

In a first time, we are going to study estimators of $f$ in every unidimensionnal model $\{\alpha f_{k}(.),\alpha\in\mathds{R}\}$
(as done in \cite{iterative}).
Usually these models are too small and the obtained estimators do not have good properties.
We then propose an iterative method that selects and aggregate such estimators in order to build a suitable estimator of $f$
(section \ref{em}).

In section \ref{rc} we study the rate of convergence of the obtained estimator in a particular case.

In section \ref{bb} we propose several improvements and generalizations of the method.

Finally, in section \ref{si} we make some simulations in order to compare the practical performances
of our estimator with other ones.

\section{Estimation method}

\label{em}

\subsection{Hypothesis}

In this section we will use a particular hypothesis about $f$ and/or the basis functions $f_{k},k\in\{1,...,m\}$.

\begin{dfn}
We will say that $f$ and $(f_{1},...,f_{m})$ satisfies the conditions $\mathcal{H}(p)$ for $1<p < +\infty$
if, for:
$$ \frac{1}{p}+\frac{1}{q} =1, $$
there exists some $(c,c_{1},...,c_{m})\in\left(\mathds{R}_{+}^{*}\right)^{m+1}$ (known to the statistician) such that:
\begin{align*}
\forall k\in\{1,...,m\}, \quad
\left(\int_{\mathcal{X}}\left|f_{k}\right|^{2p}\lambda(dx)\right)^{\frac{1}{p}} & \leq c_{k} \int_{\mathcal{X}}\left|f_{k}\right|^{2}\lambda(dx) \\
\text{and} \quad
\left(\int_{\mathcal{X}}\left|f\right|^{q}\lambda(dx)\right)^{\frac{1}{p}} & \leq c \int_{\mathcal{X}}\left|f\right|\lambda(dx)
\quad \left( = c \right)
.
\end{align*}
For $p=1$ the condition $\mathcal{H}(1)$ is: $f$ is bounded by a (known) constant $c$ and we put $c_{1}=...=c_{k}=1$.
For $p=+\infty$ the condition $\mathcal{H}(+\infty)$ is just that every $|f_{k}|$ is bounded by
$$\sqrt{c_{k}\int_{\mathcal{X}}f_{k}(x)^{2}\lambda(dx)} $$
where $c_{k}$ is known,
and we put $c=1$.
In any case, we put, for any $k$:
$$ C_{k} = c_{k} c .$$
\end{dfn}

\begin{dfn}
We put, for any $k\in\{1,...,m\}$:
$$ D_{k} = \int_{\mathcal{X}}\left|f_{k}\right|^{2}\lambda(dx) = d^{2}(f_{k},0) = \|f_{k}\|^{2} .$$
\end{dfn}

\subsection{Unidimensionnal models}

Let us choose $k\in\{1,...,m\}$ and consider the unidimensionnal model $\mathcal{M}_{k}=\{\alpha f_{k}(.),\alpha\in\mathds{R} \}$.
Remark that the orthogonal projection (denoted by $\Pi_{\mathcal{M}_{k}}$) of $f$ on $\mathcal{M}_{k}$ is known, it is namely:
$$ \Pi_{\mathcal{M}_{k}} f (.) = \overline{\alpha}_{k} f_{k}(.) $$
where:
$$ \overline{\alpha}_{k} = \arg\min_{\alpha\in\mathds{R}} d^{2} (\alpha f_{k}, f)
   = \frac{\int_{\mathcal{X}}f_{k}(x)f(x)\lambda(dx)}{\int_{\mathcal{X}}f_{k}(x)^{2}\lambda(dx)}
   = \frac{\int_{\mathcal{X}}f_{k}(x)f(x)\lambda(dx)}{D_{k}}
. $$

A natural estimator of this coefficient is:
$$ \hat{\alpha}_{k} = \frac{\frac{1}{N}\sum_{i=1}^{N} f_{k}(X_{i}) }{\int_{\mathcal{X}}f_{k}(x)^{2}\lambda(dx)},$$
because we expect to have, by the law of large numbers:
$$ \frac{1}{N}\sum_{i=1}^{N} f_{k}(X_{i}) \xrightarrow[N\rightarrow \infty]{a.s.} 
   P\left[f_{k}(X) \right] = \int_{\mathcal{X}}f_{k}(x)f(x)\lambda(dx) .$$

Actually, we can formulate a more precise result.

\begin{thm}
\label{firstbound}
Let us assume that condition $\mathcal{H}(p)$ holds for some $p\in\left[1,+\infty\right]$. Then
for any $\varepsilon>0$ we have:
\begin{multline*}
P^{\otimes N}
\Biggl\{
\forall k\in\{1,...,m\},
d^{2}(\hat{\alpha}_{k}f_{k}, \overline{\alpha}_{k} f_{k})
\\
\leq
\left\{\frac{4\left[1+\log\frac{2m}{\varepsilon}\right]}{N}\right\}
\left[\frac{\frac{1}{N}\sum_{i=1}^{N} f_{k}(X_{i})^{2} }{D_{k}} + C_{k} 
        \right]
\Biggr\}
\geq 1-\varepsilon.
\end{multline*}
\end{thm}

The proof is given at the end of the section.

\subsection{The selection algorithm}

Until the end of this section we assume that $\mathcal{H}(p)$ is satisfied for some $1\leq p\leq + \infty$.

Let $\beta(\varepsilon,k)$ denote the upper bound for the model $k$ in theorem  \ref{firstbound}:
$$ \forall \varepsilon>0, \forall k\in\{1,...,m\}: \quad
    \beta(\varepsilon,k) = \left\{\frac{4\left[1+\log\frac{2m}{\varepsilon}\right]}{N}\right\}
\left[\frac{\frac{1}{N}\sum_{i=1}^{N} f_{k}(X_{i})^{2} }{D_{k}} + C_{k}
        \right]
.$$
Let us put:
$$ \mathcal{CR}_{k,\varepsilon} = \biggl\{g\in\mathcal{L}^{2},
d^{2}(\hat{\alpha}_{k}f_{k},\Pi_{\mathcal{M}_{k}}g) \leq \beta(\varepsilon,k) \biggr\} .$$
Then theorem \ref{firstbound} implies the following result.

\begin{cor} \label{cor1}
For any $\varepsilon>0$ we have:
$$
P^{\otimes N}
\biggl\{
\forall k\in\{1,...,m\},
f\in\mathcal{CR}_{k,\varepsilon}
\biggr\}
\geq 1-\varepsilon.
$$
\end{cor}

So for any $k$, $\mathcal{CR}_{k,\varepsilon}$ is a confidence region at level $k$ for $f$.
Moreover, $\mathcal{CR}_{k,\varepsilon}$ being convex we have the following corollary.

\begin{cor} \label{cor2}
For any $\varepsilon>0$ we have:
$$
P^{\otimes N}
\biggl\{
\forall k\in\{1,...,m\}, \forall g\in \mathcal{L}^{2},
d^{2}(\Pi_{\mathcal{CR}_{k,\varepsilon}}g,f) \leq d^{2}(g,f)
\biggr\}
\geq 1-\varepsilon.
$$
\end{cor}

It just means that for any $g$, $\Pi_{\mathcal{M}_{k}}g$ is a better estimator than $g$.

So we propose the following algorithm (generic form):
\begin{itemize}
\item we choose $\varepsilon$ and start with $g_{0}=0$;
\item at each step $n$, we choose a model $\mathcal{M}_{k(n)}$ where $k(n)\in\{1,...,m\}$ can be chosen
      on each way we want (it can of course depend on the data) and take:
      $$ g_{n+1} = \Pi_{\mathcal{CR}_{k(n),\varepsilon}}g_{n} ;$$
\item we choose a stopping time $n_{s}$ on each way we want and take:
      $$\hat{f} = g_{n_{s}}.$$
\end{itemize}

So corollary \ref{cor2} implies that:
$$ P^{\otimes N}
\biggl\{
d^{2}(\hat{f},f) = d^{2}(g_{n_{s}},f) \leq ... \leq d^{2}(g_{0},f)= d^{2}(0,f)
\biggr\}
\geq 1-\varepsilon.
$$

Actually, a more accurate version of corollary \ref{cor2} can give an idea of the way to choose $k(n)$ in the algorithm.
Let us use corollary \ref{cor1} and remember the fact that each $\mathcal{CR}_{k,\varepsilon}$ is convex.

\begin{cor} \label{cor3}
For any $\varepsilon>0$ we have:
\begin{multline*}
P^{\otimes N}
\biggl\{
\forall k\in\{1,...,m\}, \forall g\in \mathcal{L}^{2},
d^{2}(\Pi_{\mathcal{CR}_{k,\varepsilon}}g,f) \leq d^{2}(g,f) - d^{2}(\Pi_{\mathcal{CR}_{k,\varepsilon}}g,g)
\biggr\}
\\
\geq 1-\varepsilon.
\end{multline*}
\end{cor}

So we propose the following version of our previous algorithm (this is not necessarily the better choice!):
\begin{itemize}
\item we choose $\varepsilon$ and $0<\kappa\leq 1/N$ and start with $g_{0}=0$;
\item at each step $n$, we take:
      $$k(n) = \arg\max_{k\in\{1,...,m\}} d^{2}(\Pi_{\mathcal{CR}_{k,\varepsilon}}g_{n},g_{n}) $$
      and:
      $$ g_{n+1} = \Pi_{\mathcal{CR}_{k(n),\varepsilon}}g_{n} ;$$
\item we take:
      $$n_{s}=\inf \left\{n\in\mathds{N}:\quad d^{2}(g_{n},g_{n-1}) \leq \kappa \right\} $$
      and:
      $$\hat{f} = g_{n_{s}}.$$
\end{itemize}

So corollary \ref{cor3} implies that:
$$ P^{\otimes N}
\biggl\{
d^{2}(\hat{f},f) \leq d^{2}(0,f) - \sum_{n=0}^{n_{s}-1}d^{2}(g_{n},g_{n+1})
\biggr\}
\geq 1-\varepsilon.
$$

\subsection{Remarks on the intersection of the confidence regions} \label{ag}

Actually, corollary \ref{cor1} could motivate another method. Note that:
$$ \forall k\in\{1,...,m\}, f\in\mathcal{RC}_{k,\varepsilon} \Leftrightarrow f\in\bigcap_{k=1}^{m} \mathcal{RC}_{k,\varepsilon} .$$
Let us put, for any $I\subset\{1,...,m\}$:
$$ \mathcal{RC}_{I,\varepsilon} = \bigcap_{k\in I} \mathcal{RC}_{k,\varepsilon} ,$$
and:
$$ \hat{f}_{I} = \Pi_{\mathcal{RC}_{I,\varepsilon}} 0 .$$

Then $ \mathcal{RC}_{I,\varepsilon}$ is still a convex region that contains $f$ and is a subset of every $\mathcal{RC}_{k,\varepsilon}$ for $k\in I$.
So we have the following result.

\begin{cor} \label{cor4}
For any $\varepsilon>0$:
\begin{multline*}
P^{\otimes N}
\biggl\{
\forall I \subset \{1,...,m\}, \forall k\in I,
d(\hat{f}_{\{1,...,m\}},f) \leq d(\hat{f}_{I},f) \leq d(\Pi_{\mathcal{RC}_{k,\varepsilon}} 0, f)
\biggr\}
\\
\geq 1-\varepsilon.
\end{multline*}
\end{cor}

In the case where we are interested in "model selection type aggregation" of estimators, note that, with probability at least $1-\varepsilon$:
$$
d(\Pi_{\mathcal{RC}_{k,\varepsilon}} 0, f)
\leq
d(\Pi_{\mathcal{RC}_{k,\varepsilon}} 0, \overline{\alpha}_{k}f_{k})
+
d(\overline{\alpha}_{k}f_{k},f)
\leq
\beta(\varepsilon,k) + d(f_{k},f).
$$
So we have the following result.
\begin{cor} \label{cor5}
For any $\varepsilon>0$:
$$
P^{\otimes N}
\biggl\{
d(\hat{f}_{\{1,...,m\}},f) \leq \inf_{k\in\{1,...,m\}} \left[d(f_{k},f) +\beta(\varepsilon,k) \right]
\biggr\}
\geq 1-\varepsilon.
$$
\end{cor}

The estimator $\hat{f}_{1,...,m}$ can be reached by solving the following optimization problem:

\begin{align*}
& \min_{g\in\mathcal{L}^{2}} \|g\|^{2},
\\
& s. t. \quad \forall k\in\{1,...,m\}:
\\
&
\left\{
\begin{array}{l}
\left<g-\hat{\alpha}_{k}f_{k},f_{k}\right> - \sqrt{D_{k} \beta(\varepsilon,k)} \leq 0,
\\
- \left<g-\hat{\alpha}_{k}f_{k},f_{k}\right> - \sqrt{D_{k} \beta(\varepsilon,k)}  \leq 0.
\end{array}
\right.
\end{align*}

The problem can be solved in dual form:
$$
\max_{\gamma\in\mathds{R}^{m}}
\left[
-\sum_{i=1}^{m} \sum_{k=1}^{m} \gamma_{i} \gamma_{k} \left<f_{i},f_{k}\right>
+ 2\sum_{k=1}^{m} \gamma_{k} \hat{\alpha}_{k}\|f_{k}\|^{2}
\\
-2\sum_{k=1}^{m} \left|\gamma_{k}\right| \sqrt{D_{k} \beta(\varepsilon,k)}
\right].
$$
with solution $\gamma^{*}=(\gamma_{1}^{*},...,\gamma_{m}^{*})$ and:
$$ \hat{f}_{\{1,...,m\}} = \sum_{k=1}^{m} \gamma_{k}^{*} f_{k} . $$

As:
$$ -\sum_{i=1}^{m} \sum_{k=1}^{m} \gamma_{i}^{*} \gamma_{k}^{*} \left<f_{i},f_{k}\right>
      = \left\|f^{*}\right\|^{2} $$
and:
$$ 2\sum_{k=1}^{m} \gamma^{*}_{k} \hat{\alpha}_{k}\|f_{k}\|^{2}
    = 2 \sum_{k=1}^{m} \gamma^{*}_{k} \frac{\frac{1}{N}\sum_{i=1}^{N}f_{k}(X_{i})}{\|f_{k}\|^{2}} \|f_{k}\|^{2}
    = \frac{2}{N} \sum_{i=1}^{N} f^{*}(X_{i})
$$
we can see this as a penalized maximization of the likelihood.

We can note that it is easier and more computationnaly efficient to project successively on every region $\mathcal{RC}(k,\varepsilon)$
than to project once on $\mathcal{RC}(\{1,...,m\},\varepsilon)$.

\subsection{An example: the histogram}

Let us assume that $\lambda$ is a finite measure and let $A_{1},...,A_{m}$ be a partition of $\mathcal{X}$.
We put, for any $k\in\{1,...,m\}$:
$$ f_{k}(.) = \mathds{1}_{A_{k}}(.) .$$

Remark that:
$$ D_{k} = \int_{\mathcal{X}}f_{k}(x)^{2}\lambda(dx) = \lambda\left(A_{k}\right) ,$$
and that condition $\mathcal{H}(+\infty)$ is satisfied with constants:
$$ c_{k} = \frac{1}{\lambda\left(A_{k}\right)} $$
and (as we have the convention $c=1$ in this case) $C_{k} = c_{k} c = c_{k} $.

In this context we have:
\begin{align*}
\overline{\alpha}_{k} & = \frac{P\left(X\in A_{k}\right)}{\lambda\left(A_{k}\right)},  \\
\hat{\alpha}_{k} & = \frac{\frac{1}{N} \sum_{i=1}^{N}\mathds{1}_{A_{k}}\left(X_{i}\right) }{\lambda\left(A_{k}\right)}, \\
\beta(\varepsilon,k) & = \left\{\frac{4\left[1+\log\frac{2m}{\varepsilon}\right]}{N \lambda\left(A_{k}\right)}\right\}
\left[\frac{1}{N}\sum_{i=1}^{N} f_{k}(X_{i})^{2} + 1
        \right].
\end{align*}

Finally, note that all the confidence regions $ \mathcal{CR}_{k,\varepsilon}$ are all orthogonal in this case. So the order
of projection does not affect the obtained estimator here, and we can take:
$$ \hat{f}  = \Pi_{\mathcal{CR}_{m,\varepsilon}} ... \Pi_{\mathcal{CR}_{1,\varepsilon}} 0 $$
(and note that $\hat{f}=\hat{f}_{\{1,...,m\}}$ here, following the notations of subsection \ref{ag}).
We have:
$$
\hat{f}(x) = \sum_{k=1}^{m} \left(\hat{\alpha}_{k} - \sqrt{\lambda\left(A_{k}\right)\beta(\varepsilon,k)} \right)_{+} f_{k}(x)
$$
where, for any $y\in\mathds{R}$:
$$
(y)_{+} = \left\{
\begin{array}{l}
y \quad \text{if} \quad y \geq 0
\\
\\
0 \quad \text{otherwise.}
\end{array}
\right.
$$

In this case corollary \ref{cor3} becomes:
$$ P^{\otimes N}
\biggl\{
d^{2}(\hat{f},f) \leq d^{2}(0,f) -
\sum_{k=1}^{m} \left(\hat{\alpha}_{k} - \sqrt{\lambda\left(A_{k}\right)\beta(\varepsilon,k)} \right)_{+} ^{2}\lambda\left(A_{k}\right)
\biggr\}
\geq 1-\varepsilon.
$$

\subsection{Proof of the theorem}

Before giving the proof, let us state two lemmas that we will use in the proof.
The first one is a variant of a lemma by Catoni \cite{Classif}, the second one
is due to Panchenko \cite{Panchenko}.

\begin{lemma}
\label{catonilemma}
Let $(T_{1},...,T_{2N})$ be a random vector taking values in $\mathds{R}^{2N}$ distributed according to
a distribution $\mathcal{P}^{\otimes 2N}$.
For any $\eta\in\mathds{R}$,
for any measurable function $\lambda:\mathds{R}^{2N}\rightarrow\mathds{R}_{+}^{*}$
that is exchangeable with respect to its $2\times 2N$ arguments:
$$
\mathcal{P}^{\otimes 2N} \exp\Biggl(\frac{\lambda}{N}\sum_{i=1}^{N}\Bigl\{ T_{i+N}-T_{i} \Bigr\}
  - \frac{\lambda^{2}}{N^{2}} \sum_{i=1}^{2N} T_{i}^{2}
  - \eta \Biggr) \leq \exp \left( - \eta \right)
$$
and the reverse inequality:
$$
\mathcal{P}^{\otimes 2N} \exp\Biggl(\frac{\lambda}{N}\sum_{i=1}^{N}\Bigl\{ T_{i} - T_{i+1} \Bigr\}
  - \frac{\lambda^{2}}{N^{2}} \sum_{i=1}^{2N} T_{i}^{2}
  - \eta \Biggr) \leq \exp \left( - \eta \right)
,
$$
where we write:
\begin{align*}
\eta & =\eta\left(T_{1},...,T_{2N}\right) \\
\lambda & =\lambda\left(T_{1},...,T_{2N}\right)
\end{align*}
for short.
\end{lemma}

\begin{proof}[Proof of lemma \ref{catonilemma}]
In order to prove the first inequality, we write:
\begin{multline*}
\mathcal{P}^{\otimes 2N}
\exp\Biggl(\frac{\lambda}{N}\sum_{i=1}^{N}\Bigl\{ T_{i+N}-T_{i} \Bigr\}
  - \frac{\lambda^{2}}{N^{2}} \sum_{i=1}^{2N} T_{i}^{2}
  - \eta \Biggr)
\\
=
\mathcal{P}^{\otimes 2N} \exp\Biggl(\sum_{i=1}^{N} \log \cosh \left\{ \frac{\lambda}{N} \left(T_{i+1}-T_{i}\right) \right\}
  - \frac{\lambda^{2}}{N^{2}} \sum_{i=1}^{2N} T_{i}^{2}
  - \eta \Biggr).
\end{multline*}
We now use the inequality:
$$ \forall x\in\mathds{R}, \log\cosh x \leq \frac{x^{2}}{2} .$$
We obtain:
\begin{equation*}
\log\cosh \left\{\frac{\lambda}{N}\left(T_{i+1}-T_{i}\right)  \right\}
\leq \frac{\lambda^{2}}{2N^{2}} \left(T_{i+1}-T_{i}\right)^{2}
\leq \frac{\lambda^{2}}{N^{2}} \left(T_{i+1}^{2}+T_{i}^{2}\right)
.
\end{equation*}
The proof for the reverse inequality is exactly the same.
\end{proof}

\begin{lemma}[Panchenko \cite{Panchenko}, corollary 1]
\label{panlemma}
Let us assume that we have i.i.d. variables $T_{1},...,T_{N}$ (with distribution $\mathcal{P}$ and values in $\mathds{R}$)
and an independant copy $T'=(T_{N+1},...,T_{2N})$ of 
$T=(T_{1},...,T_{N})$. Let $\xi_{j}(T,T')$ for $j\in\{1,2,3\}$ be three measurables functions taking values in $\mathds{R}$, and
$\xi_{3}\geq 0$. Let us assume that we know two constants $A\geq 1$ and $a>0$ such that, for any $u>0$:
$$ P^{\otimes 2N} \left[\xi_{1}(T,T')\geq \xi_{2}(T,T')+\sqrt{\xi_{3}(T,T')u}\right] \leq A\exp(-au) .$$
Then, for any $u>0$:
\begin{multline*}
P^{\otimes 2N} \Biggl\{P^{\otimes 2N}\left[\xi_{1}(T,T')|T\right]
\\
         \geq P^{\otimes 2N}\left[\xi_{2}(T,T')|T\right]+
          \sqrt{P^{\otimes 2N}\left[\xi_{3}(T,T')|T\right]u}\Biggr\} \leq A\exp(1-au) .
\end{multline*}
\end{lemma}

The proof of this lemma can be found in \cite{Panchenko}.
We can now give the proof of theorem \ref{firstbound}.

\begin{proof}[Proof of theorem \ref{firstbound}]
Let $(X_{N+1},...,X_{2N})$ be an independant copy of our sample $(X_{1},...,X_{N})$.
Let us choose $k\in\{1,...,m\}$.
Let us apply lemma \ref{catonilemma} with $\mathcal{P}=P$ and, for any $i\in\{1,...,2N\}$:
$$ T_{i} = f_{k}(X_{i}).$$
We obtain, for any
measurable function $\eta_{k}\in\mathds{R}$,
for any measurable function $\lambda_{k}:\mathds{R}^{2N}\rightarrow\mathds{R}_{+}^{*}$
that is exchangeable with respect to its $2\times 2N$ arguments:
$$
P^{\otimes 2N} \exp\Biggl(\frac{\lambda_{k}}{N}\sum_{i=1}^{N}\Bigl\{ f_{k}(X_{i+N})-f_{k}(X_{i}) \Bigr\}
 - \frac{\lambda_{k}^{2}}{N^{2}} \sum_{i=1}^{2N} f_{k}(X_{i})^{2}
  - \eta_{k} \Biggr) \leq \exp \left( - \eta_{k} \right)
$$
and the reverse inequality:
$$
P^{\otimes 2N} \exp\Biggl(\frac{\lambda_{k}}{N}\sum_{i=1}^{N}\Bigl\{ f_{k}(X_{i}) - f_{k}(X_{i+N}) \Bigr\}
  - \frac{\lambda_{k}^{2}}{N^{2}} \sum_{i=1}^{2N} f_{k}(X_{i})^{2}
  - \eta_{k} \Biggr) \leq \exp \left( - \eta_{k} \right)
$$
as wall. This implies that:
$$
P^{\otimes 2N}
\Biggl[
\frac{1}{N}\sum_{i=1}^{N}\Bigl\{ f_{k}(X_{i}) - f_{k}(X_{i+N}) \Bigr\}
\leq
\frac{\lambda_{k}}{N^{2}} \sum_{i=1}^{2N} f_{k}(X_{i})^{2}
  + \frac{\eta_{k}}{\lambda_{k}}
\Biggr]
\leq
\exp \left( - \eta_{k} \right)
$$
and:
$$
P^{\otimes 2N}
\Biggl[
\frac{1}{N}\sum_{i=1}^{N}\Bigl\{ f_{k}(X_{i+N}) - f_{k}(X_{i}) \Bigr\}
\leq
\frac{\lambda_{k}}{N^{2}} \sum_{i=1}^{2N} f_{k}(X_{i})^{2}
  + \frac{\eta_{k}}{\lambda_{k}}
\Biggr]
\leq
\exp \left( - \eta_{k} \right).
$$
Let us choose:
$$\lambda_{k} = \sqrt{\frac{N\eta_{k}}{\frac{1}{N}\sum_{i=1}^{2N}f_{k}(X_{i})^{2}}}$$
in both inequalities, we obtain for the first one:
$$
P^{\otimes 2N}
\Biggl[
\frac{1}{N}\sum_{i=1}^{N}\Bigl\{ f_{k}(X_{i}) - f_{k}(X_{i+N}) \Bigr\}
\geq
2\sqrt{\frac{\eta_{k}\frac{1}{N}\sum_{i=1}^{2N}f_{k}(X_{i})^{2}}{N}}
\Biggr]
\leq
\exp \left( - \eta_{k} \right).
$$
We now apply lemma \ref{panlemma} with the same $T_{i}=f_{k}(X_{i})$, $\eta_{k}=u$,
$A=1$, $a=1$, $\xi_{2}=0$,
$$ \xi_{1} = \frac{1}{N}\sum_{i=1}^{N}\Bigl\{ f_{k}(X_{i}) - f_{k}(X_{i+N}) \Bigr\} \quad \text{and}$$
$$ \xi_{3} = \frac{4\frac{1}{N}\sum_{i=1}^{2N}f_{k}(X_{i})^{2}}{N}.$$
We obtain:
\begin{multline*}
P^{\otimes N}
\Biggl[
\frac{1}{N}\sum_{i=1}^{N}f_{k}(X_{i})-P\left[f_{k}(X)\right]
\geq
2\sqrt{\frac{\eta_{k}\left\{\frac{1}{N}\sum_{i=1}^{N}f_{k}(X_{i})^{2}+P\left[f_{k}(X)^{2}\right]\right\}}{N}}
\Biggr]
\\
=
P^{\otimes 2N}
\Biggl[
\frac{1}{N}\sum_{i=1}^{N}f_{k}(X_{i})-P\left[f_{k}(X)\right]
\geq
2\sqrt{\frac{\eta_{k}\left\{\frac{1}{N}\sum_{i=1}^{N}f_{k}(X_{i})^{2}+P\left[f_{k}(X)^{2}\right]\right\}}{N}}
\Biggr]
\\
\leq
\exp \left( 1 - \eta_{k} \right).
\end{multline*}
Remark that:
$$
P\left[f_{k}(X)^{2}\right]
    = \int_{\mathcal{X}} f_{k}(x)^{2} f(x)\lambda(dx).
$$
So, using condition $\mathcal{H}(p)$ and Hölder's inequality we have:
\begin{multline*}
P\left[f_{k}(X)^{2}\right] \leq \left(\int_{\mathcal{X}} \left|f_{k}(x)\right|^{2p} \lambda(dx)\right)^{\frac{1}{p}}
                                \left(\int_{\mathcal{X}} f(x)^{q} \lambda(dx)\right)^{\frac{1}{q}}
\\
\leq
\left(c_{k} \int_{\mathcal{X}} f_{k}(x)^{2} \lambda(dx)\right)\left( c\int_{\mathcal{X}} f(x)\lambda(dx)\right)
\\
= \left(c_{k}c\right) \int_{\mathcal{X}} f_{k}(x)^{2} \lambda(dx) = C_{k} D_{k}.
\end{multline*}
Now, let us combine this inequality with the reverse one by a union bound argument, we have:
\begin{multline*}
P^{\otimes N}
\Biggl[
\left|\frac{1}{N}\sum_{i=1}^{N}f_{k}(X_{i})-P\left[f_{k}(X)\right]\right|
\\
\geq
2\sqrt{\frac{\eta_{k}\left\{\frac{1}{N}\sum_{i=1}^{N}f_{k}(X_{i})^{2}+C_{k}D_{k}\right\}}{N}}
\Biggr]
\leq
2 \exp \left( 1 - \eta_{k} \right).
\end{multline*}
We now make a union bound on $k\in\{1,...,m\}$ and put:
$$ \eta_{k} = 1+\log\frac{2m}{\varepsilon}. $$
We obtain:
\begin{multline*}
P^{\otimes N}
\Biggl[
\forall k\in\{1,...,m\}, \quad
\left|\frac{1}{N}\sum_{i=1}^{N}f_{k}(X_{i})-P\left[f_{k}(X)\right]\right|
\\
\leq
2\sqrt{\frac{\left(1+\log\frac{2m}{\varepsilon}\right)\left\{\frac{1}{N}\sum_{i=1}^{N}f_{k}(X_{i})^{2}+C_{k}D_{k}\right\}}{N}}
\Biggr]
\geq
1-\varepsilon.
\end{multline*}
We end the proof by noting that:
$$
d^{2}(\hat{\alpha}_{k}f_{k}, \overline{\alpha}_{k} f_{k})
= \frac{
\left[\frac{1}{N}\sum_{i=1}^{N}f_{k}(X_{i})-P\left[f_{k}(X)\right]\right]^{2}
} { \int_{\mathcal{X}}f_{k}(x)^{2}\lambda(dx) }.
$$
\end{proof}

\section{Some examples with rates of convergence}

\label{rc}

\subsection{General remarks when $(f_{k})_{k}$ is an orthonormal family and condition $\mathcal{H}(1)$ is satisfied}

\label{rcgr}

In subsections \ref{rcgr}, \ref{sob} and \ref{bes}, we study the rate of convergence of our estimator in the special case where
$(f_{k})_{k\in\mathds{N}^{*}}$ is an orthonormal basis of $\mathcal{L}^{2}$, so we have:
$$ D_{k} = \int_{\mathcal{X}} f_{k}(x)^{2}\lambda(dx) = 1 $$
and:
$$ \int_{\mathcal{X}} f_{k}(x)f_{k'}(x)\lambda(dx) = 0 $$
if $k\neq k'$.

We also assume that condition $\mathcal{H}(1)$ is satisfied: $\forall x\in\mathcal{X}, f(x) \leq c$, remember that in this case
we have taken $c_{k}=1$ and so $C_{k} = c$, so:
$$ \beta(\varepsilon,k) = \left\{\frac{4\left[1+\log\frac{2m}{\varepsilon}\right]}{N}\right\}
\left[\frac{1}{N}\sum_{i=1}^{N} f_{k}(X_{i})^{2} + c
        \right]. $$

Note that in this case all the order of application of the projections $\Pi_{\mathcal{RC}_{k,\varepsilon}}$ does not matter
because these projections works on orthogonal directions. So we can define, once $m$ is chosen:
$$ \hat{f} = \Pi_{\mathcal{RC}_{m,\varepsilon}} ... \Pi_{\mathcal{RC}_{1,\varepsilon}} 0
     = \Pi_{\mathcal{RC}_{\{1,...,m\},\varepsilon}} 0 = \hat{f}_{\{1,...,m\}} $$
(following the notations of subsection \ref{ag}).
Note that:
$$ \hat{f}(x) = \sum_{k=1}^{m} sign(\hat{\alpha}_{k})\left(\left|\hat{\alpha}_{k}\right|-\sqrt{\beta(\varepsilon,k)}\right)_{+}f_{k}(x) $$
where $sign(x)$ is the sign of $x$ (namely $+1$ if $x>0$ and $-1$ otherwise), and
so $\hat{f}$ is a soft-thresholded estimator.
Let us also make the following remark. As for any $x$, $f(x)\leq c$, we have:
$$ d^{2}(f,0) \leq c .$$
So the region:
$$ \mathcal{B} = \left\{g\in\mathcal{L}^{2}: \forall k\in\mathds{N}^{*}, \int_{\mathcal{X}}g(x)f_{k}(x)\lambda(dx)\leq \sqrt{c}\right\}$$
is convex, and contains $f$. So the projection on $\mathcal{B}$, $\Pi_{\mathcal{B}}$ can only improve $\hat{f}$.
We put:
$$ \tilde{f} = \Pi_{\mathcal{B}}\hat{f}.$$
Note that this transormation is needed to obtain the following theorem, but does not have practical incidence in general.
Actually:
$$ \tilde{f}(x)
= \sum_{k=1}^{m} sign(\hat{\alpha}_{k}) \left\{ \left(\left|\hat{\alpha}_{k}\right|-\sqrt{\beta(\varepsilon,k)}\right)_{+}
      \wedge\sqrt{c} \right\}f_{k}(x) .$$

\subsection{Rate of convergence in Sobolev spaces}

\label{sob}

It is well known that if $f$ has regularity $\beta$ (known by the statistician) then we have the choice
$$m=N^{\frac{1}{2\beta+1}}$$
and a standard estimation of coefficients leads to the optimal rate of convergence:
$$ N^{\frac{-2\beta}{2\beta+1}}.$$

Here, we assume that we don't know $\beta$, and we show that taking $m=N$ leads
to the rate of convergence:
$$ N^{\frac{-2\beta}{2\beta+1}} \log N $$
namely the optimal rate of convergence up to a $\log N$ factor.

\begin{thm} \label{rate1}
Let us assume that $(f_{k})_{k\in\mathds{N}^{*}}$ is an orthonormal basis of $\mathcal{L}^{2}$.
Let us put:
$$\overline{f}_{m} =\arg\min_{g\in Span(f_{1},...,f_{m})} d^{2}(g,f),$$
and let us assume that $f\in\mathcal{L}^{2}$ satisfies condition $\mathcal{H}(1)$ and
is such that there are unknown constants $D>0$ and $\beta\geq 1$ such that:
$$ d^{2}(\overline{f}_{m},f) \leq D m^{-2\beta} .$$
Let us choose $m=N$ and $\varepsilon=N^{-2}$ in the definition of $\tilde{f}$.
Then we have, for any $N\geq 2$:
$$P^{\otimes N} d^{2}(\tilde{f},f) \leq D'(c,D) \left(\frac{\log N}{N}\right)^{\frac{2\beta}{2\beta+1}} .$$
\end{thm}

Here again, the proof of the theorems are given at the end of the section.
Let us just remark that,
in the case where $\mathcal{X}=[0,1]$, $\lambda$ is the Lebesgue measure, and $(f_{k})_{k\in\mathds{N}^{*}}$ is the trigonometric basis,
the condition:
$$ d^{2}(\overline{f}_{m},f) \leq D m^{-2\beta} $$
is satisfied for $D=D(\beta,L)=L^{2}\pi^{-2\beta}$ as soon as $f\in W(\beta,L)$ where $W(\beta,L)$ is the Sobolev class:
$$ \left\{f\in\mathcal{L}^{2}: f^{(\beta-1)} \text{ is absolutely continuous and }
\int_{0}^{1}f^{(\beta)}(x)^{2}\lambda(dx) \leq L^{2}
\right\}, $$
see Tsybakov \cite{nonpara} for example.
The minimax rate of convergence in $W(\beta,L)$ is $N^{-\frac{2\beta}{2\beta+1}}$,
so we can see that our estimator reaches the best rate of convergence up to a $\log N$ factor
with an unknown $\beta$.

\subsection{Rate of convergence in Besov spaces}

\label{bes}

We here extend the previous result to the case of a Besov space $B_{s,p,q}$.
Note that we have, for any $L\geq 0$ and $\beta\geq 0$:
$$ W(\beta,L) \subset B_{\beta,2,2} $$
so this result is really an extension of the previous one (see Härdle, Kerkyacharian, Picard and Tsybakov \cite{Ondel},
or Donoho, Johnstone, Kerkyacharian and Picard \cite{Ondel2}).
We define the Besov space:
\begin{multline*}
B_{s,p,q} = \Biggl\{g:[0,1]\rightarrow \mathds{R}, \quad
g(.)=\alpha\phi(.) + \sum_{j=0}^{\infty}\sum_{k=1}^{2^{j}}\beta_{j,k}\psi_{j,k}(.), \quad
\\
\sum_{j=0}^{\infty}2^{jq\left(s-\frac{1}{2}-\frac{1}{p}\right)} \left[\sum_{k=1}^{2^{j}} \left|\beta_{j,k}\right|^{p}\right]^{\frac{q}{p}}
= \|g\|_{s,p,q}^{q} < + \infty
\Biggr\},
\end{multline*}
with obvious changes for $p=+\infty$ or $q=+\infty$.
We also define the weak Besov space:
\begin{multline*}
W_{\rho,\pi} = \Biggl\{g:[0,1]\rightarrow \mathds{R}, \quad
g(.)=\alpha\phi(.) + \sum_{j=0}^{\infty}\sum_{k=1}^{2^{j}}\beta_{j,k}\psi_{j,k}(.), \quad
\\
\sup_{\lambda>0} \lambda^{\rho} \sum_{j=0}^{\infty} 2^{j\left(\frac{\pi}{2}-1\right)} \sum_{k=1}^{2^{j}}
\mathds{1}_{\left\{|\beta_{j,k}|>\lambda\right\}}
<+\infty
\Biggr\}
\\
=\Biggl\{g:[0,1]\rightarrow \mathds{R}, \quad
g(.)=\alpha\phi(.) + \sum_{j=0}^{\infty}\sum_{k=1}^{2^{j}}\beta_{j,k}\psi_{j,k}(.), \quad
\\
\sup_{\lambda>0} \lambda^{\pi-\rho} \sum_{j=0}^{\infty} 2^{j\left(\frac{\pi}{2}-1\right)} \sum_{k=1}^{2^{j}}
|\beta_{j,k}|^{\pi} \mathds{1}_{\left\{|\beta_{j,k}|\leq \lambda\right\}}
<+\infty
\Biggr\},
\end{multline*}
see Cohen \cite{Cohen} for the equivalence of both definitions.
Let us remark that $B_{s,p,q}$ is a set of functions with regularity $s$ while $W_{\rho,\pi}$ is
a set of functions with regularity:
$$ s'=\frac{1}{2}\left(\frac{\pi}{\rho}-1\right). $$

\begin{thm} \label{rate2}
Let us assume that $\mathcal{X}=[0,1]$, and that
$(\psi_{j,k})_{j=0,...,+\infty, k\in\{1,...,2^{j}\}}$ is a wavelet basis, together with
a function $\phi$, satisfying the conditions given in \cite{Ondel2}
and having regularity $R$ (for example Daubechies' families), with $\phi$ and $\psi_{0,1}$ supported by $[-A,A]$.
Let us assume that $f\in B_{s,p,q}$ with $R+1\geq s>\frac{1}{p}$, $1\leq q\leq \infty$, $2\leq p\leq +\infty$,
or that $f\in B_{s,p,q}\cap W_{\frac{2}{2s+1},2}$ with $R+1\geq s>\frac{1}{p}$, $1\leq p\leq +\infty$,
with unknown constants $s$, $p$ and $q$ and that $f$ satisfies condition $\mathcal{H}(1)$ with a known constant $c$.
Let us choose:
$$ \{f_{1},...,f_{m}\} = \{\phi\} \cup \{\psi_{j,k},j=1,...,2^{\lfloor\frac{\log N}{\log 2}\rfloor},k=1,...,2^{j}\} $$
(so $\frac{N}{2} \leq m\leq N$)
and $\varepsilon=N^{-2}$ in the definition of $\tilde{f}$.
Then we have:
$$P^{\otimes N} d^{2}(\tilde{f},f) = \mathcal{O}\left( \left(\frac{\log N}{N}\right)^{\frac{2s}{2s+1}}
\right).$$
\end{thm}

Let us remark that we obtain nearly the same rate of convergence than in \cite{Ondel2}, namely
the minimax rate of convergence up to a $\log N$ factor.

\subsection{Kernel estimators}

Here, we assume that $\mathcal{X}=\mathds{R}$ and that $f$ is compactly supported, say by $[0,1]$.
We put, for any $m\in\mathds{N}$ and $k\in\{1,...,m\}$:
$$ f_{k}(x) = K\left(\frac{k}{m},x\right) $$
where $K$ is some function $\mathds{R}\times\mathds{R} \rightarrow \mathds{R}$
and we obtain some estimator that has the form of a kernel estimator:
$$ \hat{f}_{\{1,...,m\}}(x) = \sum_{k=1}^{m} \tilde{\alpha}_{k} K\left(\frac{k}{m},x\right) .$$
Moreover, is is possible to use a multiple kernel estimator. Let us choose $n\in\mathds{N}$, $h\in\mathds{N}$,
$h$ kernels $K_{1},...,K_{h}$ and put, for any $k=i+n*j\in\{1,...,m=hn\}$:
$$ f_{k}(x) = K_{j}\left(\frac{i}{n},x\right) .$$
We obtain a multiple kernel estimator:
$$ \hat{f}_{\{1,...,m\}}(x) = \sum_{i=1}^{n} \sum_{j=1}^{h} \tilde{\alpha}_{i+nj} K_{j} \left(\frac{i}{n},x\right) .$$

\subsection{Proof of the theorems}

\begin{proof}[Proof of theorem \ref{rate1}]
Let us begin the proof with a general $m$ and $\varepsilon$, the reason of the choice
$m=N$ and $\varepsilon=N^{-2}$ will become clear.
Let us also write $\mathcal{E}(\varepsilon)$ the event satisfied with probability at least $1-\varepsilon$ in theorem
\ref{firstbound}.
We have:
$$
P^{\otimes N} d^{2}(\tilde{f},f)
=
P^{\otimes N}\Biggl[\mathds{1}_{\mathcal{E}(\varepsilon)} d^{2}(\tilde{f},f) \Biggr]
+
P^{\otimes N}\Biggl[\left(1-\mathds{1}_{\mathcal{E}(\varepsilon)}\right)d^{2}(\tilde{f},f)\Biggr].
$$
For the first term we have:
$$
d^{2}(\tilde{f},f) \leq 2\int_{\mathcal{X}}f(x)^{2}\lambda(dx) + 2\int_{\mathcal{X}}\tilde{f}(x)^{2}\lambda(dx)
\leq 2c + 2mc = 2(m+1)c
$$
and so:
$$ P^{\otimes N}\Biggl[\left(1-\mathds{1}_{\mathcal{E}(\varepsilon)}\right)d^{2}(\hat{f},f)\Biggr] \leq 2 \varepsilon (m+1) c.$$
For the other term, just remark that under $\mathcal{E}(\varepsilon)$:
\begin{multline*}
d^{2}(\tilde{f},f)
=
d^{2}(\Pi_{\mathcal{B}}\Pi_{\mathcal{CR}_{m,\varepsilon}}...\Pi_{\mathcal{CR}_{1,\varepsilon}}0,f)
\\
\leq
d^{2}(\Pi_{\mathcal{CR}_{m,\varepsilon}}...\Pi_{\mathcal{CR}_{1,\varepsilon}}0,f)
\leq
d^{2}(\Pi_{\mathcal{CR}_{m',\varepsilon}}...\Pi_{\mathcal{CR}_{1,\varepsilon}}0,f)
\end{multline*}
for any $m'\leq m$, because of theorem \ref{firstbound}, more precisely of corollary \ref{cor2}.
And we have:
\begin{multline*}
d^{2}(\Pi_{\mathcal{M}_{m'}}...\Pi_{\mathcal{M}_{1}}0,f)
\\
\leq \sum_{k=1}^{m'}
\left\{\frac{4\left[1+\log\frac{2m}{\varepsilon}\right]}{N}\right\}
\left[\frac{1}{N}\sum_{i=1}^{N} f_{k}(X_{i})^{2} + c
        \right]
+ 
d^{2}(\overline{f}_{m},f).
\end{multline*}
So we have:
\begin{multline*}
P^{\otimes N}\Biggl[\mathds{1}_{\mathcal{E}(\varepsilon)} d^{2}(\tilde{f},f)\Biggr]
\leq
P^{\otimes N}\Biggl[d^{2}(\tilde{f},f)\Biggr]
\\
\leq
P^{\otimes N} \sum_{k=1}^{m'}
\left\{\frac{4\left[1+\log\frac{2m}{\varepsilon}\right]}{N}\right\}
\left[\frac{1}{N}\sum_{i=1}^{N} f_{k}(X_{i})^{2} + c
\right] + (m')^{-2\beta} D
\\
\leq
\frac{8m' c\left[1+\log\frac{2m}{\varepsilon}\right]}{N}
+
(m')^{-2\beta} D.
\end{multline*}
So finally, we obtain, for any $m'\leq m$:
$$
P^{\otimes N} d^{2}(\tilde{f},f)
\leq
\frac{8m' c\left[1+\log\frac{2m}{\varepsilon}\right]}{N}
+
(m')^{-2\beta} D
+ 2 \varepsilon (m+1) c.
$$
The choice of:
$$m'=\left(\frac{N}{\log N}\right)^{\frac{1}{2\beta+1}}$$
leads to a first term of order $N^{\frac{-2\beta}{2\beta+1}}\log \frac{m}{\varepsilon} (\log N)^{-\frac{1}{2\beta+1}}$
and a second term of order $N^{\frac{-2\beta}{2\beta+1}}(\log N)^{\frac{2\beta}{2\beta+1}}$.
The choice of $m=N$ and $\varepsilon=N^{-2}$ gives a first and second term at order:
$$ \left(\frac{\log N}{N}\right)^{\frac{2\beta}{2\beta+1}}$$
while keeping the third term at order $N^{-1}$.
This proves the theorem.
\end{proof}

\begin{proof}[Proof of theorem \ref{rate2}]
Here again let us write $\mathcal{E}(\varepsilon)$ the event satisfied with probability at least $1-\varepsilon$ in theorem
\ref{firstbound}.
We have:
$$
P^{\otimes N} d^{2}(\tilde{f},f)
=
P^{\otimes N}\Biggl[\mathds{1}_{\mathcal{E}(\varepsilon)} d^{2}(\tilde{f},f) \Biggr]
+
P^{\otimes N}\Biggl[\left(1-\mathds{1}_{\mathcal{E}(\varepsilon)}\right)d^{2}(\tilde{f},f)\Biggr].
$$
For the first term we still have:
$$
d^{2}(\tilde{f},f) \leq 2(m+1)c.
$$
For the second term, let us write the development of $f$ into our wavelet basis:
$$ f=\alpha \phi + \sum_{j=0}^{\infty}\sum_{k=1}^{2^{j}}\beta_{j,k}\psi_{j,k} ,$$
and:
$$ \hat{f}(x) = \tilde{\alpha} \phi + \sum_{j=0}^{J}\sum_{k=1}^{2^{j}}\tilde{\beta}_{j,k}\psi_{j,k} $$
the estimator $\hat{f}$.
Let us put:
$$J=\left\lfloor \frac{\log N}{\log 2} \right\rfloor.$$
For any $J'\leq J$ we have:
\begin{multline*}
d^{2}(\tilde{f},f)
=
d^{2}(\Pi_{\mathcal{B}}\Pi_{\mathcal{CR}_{m,\varepsilon}}...\Pi_{\mathcal{CR}_{1,\varepsilon}}0,f)
\leq
d^{2}(\Pi_{\mathcal{CR}_{m,\varepsilon}}...\Pi_{\mathcal{CR}_{1,\varepsilon}}0,f)
\\
=
(\tilde{\alpha} - \alpha)^{2} + \sum_{j=0}^{J}\sum_{k=1}^{2^{j}}(\tilde{\beta}_{j,k}-\beta_{j,k})^{2}
+ \sum_{j=J+1}^{\infty}\sum_{k=1}^{2^{j}}\beta_{j,k}^{2}
\\
\leq
(\tilde{\alpha} - \alpha)^{2} + \sum_{j=0}^{J'}\sum_{k=1}^{2^{j}}(\tilde{\beta}_{j,k}-\beta_{j,k})^{2}\mathds{1}(|\beta_{j,k}|\geq \kappa)
+ \sum_{j=0}^{J'}\sum_{k=1}^{2^{j}} \beta_{j,k}^{2} \mathds{1}(|\beta_{j,k}|< \kappa)
\\
+
\sum_{j=J'+1}^{\infty}\sum_{k=1}^{2^{j}}\beta_{j,k}^{2}
\end{multline*}
for any $\kappa\geq 0$, as soon as $\mathcal{E}(\varepsilon)$ is satisfied (here again we applied theorem
\ref{firstbound}).
In the case where $p\geq 2$ we can take:
$$ J'=\left\lfloor \frac{\log N^{\frac{1}{1+2s}}}{\log 2} \right\rfloor$$
and $\kappa=0$ to obtain (let $C$ be a generic constant in the whole proof):
$$
\sum_{j=J'+1}^{\infty}\sum_{k=1}^{2^{j}} \beta_{j,k}^{2}
\leq
\sum_{j=J'+1}^{\infty} \left(\sum_{k=1}^{2^{j}} \beta_{j,k}^{p}\right)^{\frac{2}{p}} 2^{j\left(1-\frac{2}{p}\right)}.
$$
As $f\in B_{s,p,q} \subset B_{s,p,\infty}$ we have:
$$ \left(\sum_{k=1}^{2^{j}} \beta_{j,k}^{p}\right)^{\frac{2}{p}} \leq C 2^{-2j\left(s+\frac{1}{2}-\frac{1}{p}\right)} $$
and so:
$$ \sum_{j=J+1}^{\infty}\sum_{k=1}^{2^{j}} \beta_{j,k}^{2} \leq C 2^{-2J's} \leq C N^{\frac{-2s}{1+2s}},$$
and:
\begin{multline*}
\sum_{j=0}^{J}\sum_{k=1}^{2^{j}}(\tilde{\beta}_{j,k}-\beta_{j,k})^{2}\mathds{1}(|\beta_{j,k}|\geq \kappa)
\leq
\frac{8c \left[1+\log\frac{2m}{\varepsilon}\right]}{N} \sum_{j=0}^{J}\sum_{k=1}^{2^{j}} 1
\\
\leq
\frac{8c \left[1+\log\frac{2m}{\varepsilon}\right]}{N} 2^{J'+1}
\leq
C \frac{8c \left[1+\log\frac{2m}{\varepsilon}\right]}{N} N^{\frac{1}{1+2s}}.
\end{multline*}
So we obtain the desired rate of convergence.
In the case where $p<2$ we let $J'=J$ and proceed as follows.
\begin{multline*}
\sum_{j=0}^{J}\sum_{k=1}^{2^{j}}(\tilde{\beta}_{j,k}-\beta_{j,k})^{2}\mathds{1}(|\beta_{j,k}|\geq \kappa)
\leq
\frac{8c \left[1+\log\frac{2m}{\varepsilon}\right]}{N} \sum_{j=0}^{J}\sum_{k=1}^{2^{j}}\mathds{1}(|\beta_{j,k}|\geq \kappa)
\\
\leq
\frac{8c \left[1+\log\frac{2m}{\varepsilon}\right]}{N} C \kappa^{-\frac{2}{2s+1}}
\end{multline*}
because $f$ is also assumed to be in the weak Besov space.
We also have:
$$
\sum_{j=0}^{J}\sum_{k=1}^{2^{j}} \beta_{j,k}^{2} \mathds{1}(|\beta_{j,k}|< \kappa)
\leq C \kappa^{2-\frac{2}{1+2s}}
.
$$
For the remainder term we use (see \cite{Ondel,Ondel2}):
$$ B_{s,p,q} \subset B_{s-\frac{1}{p}+\frac{1}{2},2,q} $$
to obtain:
$$ \sum_{j=J+1}^{\infty}\sum_{k=1}^{2^{j}} \beta_{j,k}^{2} \leq C 2^{-2J\left(s+\frac{1}{2}-\frac{1}{p}\right)}
 \leq C 2^{-J} $$
as $s>\frac{1}{p}$.
Let us remember that:
$$\frac{N}{2} \leq m = 2^{J} \leq N$$
and that $\varepsilon=N^{-2}$, and take:
$$ \kappa=\sqrt{\frac{\log N}{N}}$$
to obtain the desired rate of convergence.
\end{proof}

\section{Better bounds and generalizations}

\label{bb}

Actually, as pointed out by Catoni \cite{Classif}, the symmetrization technique used in the proof
of theorem \ref{firstbound} causes the loss of a factor $2$ in the bound because we upper bound the variance
of two samples instead of $1$.
In this section, we try to use this remark to improve our bound, using techniques already used by Catoni \cite{Cat7}.
We also give a generalization of the obtained result that allows us to use a family $(f_{1},...,f_{m})$ of functions
that is data-dependant. The technique used is due to Seeger \cite{Seeger}, and it will allows us to use kernel estimators as Support
Vector Machines.

Remark that the estimation technique described in section \ref{em}
does not necessarily require a bound on $d^{2}(\hat{\alpha}_{k}f_{k},\overline{\alpha}_{k}f_{k})$.
Actually, a simple confidence interval on $\overline{\alpha}_{k}$ is sufficient.

\subsection{An improvement of theorem \ref{firstbound} under condition $\mathcal{H}(+\infty)$}

Let us remember that $\mathcal{H}(+\infty)$ just means that every $f_{k}$ is bounded by $\sqrt{C_{k}D_{k}}$.

\begin{thm}
Under condition $\mathcal{H}(+\infty)$,
for any $\varepsilon>0$,
for any $\beta_{k,1},\beta_{k,2}$ such that:
$$  0 < \beta_{k,j} < \frac{N}{\sqrt{C_{k}D_{k}}} , \quad j\in\{1,2\}, $$
with $P^{\otimes N}$-probability at least $1-\varepsilon$,
for any $k\in\{1,...,m\}$
we have:
$$ \alpha^{\inf}_{k}(\varepsilon,\beta_{k,1}) \leq \overline{\alpha}_{k} \leq \alpha^{\sup}_{k}(\varepsilon,\beta_{k,2}) $$
with:
$$ \alpha^{\sup}_{k}(\varepsilon,\beta_{k,2}) =
\frac{
N-N\exp\left[\frac{1}{N}\sum_{i=1}^{N} \log \left(1-\frac{\beta_{k,2}}{N} f_{k}(X_{i})\right) - \frac{\log\frac{2m}{\varepsilon}}{N} \right]
}{D_{k}\beta_{k,2}}
$$
and:
$$ \alpha^{\inf}_{k}(\varepsilon,\beta_{k,1}) =
\frac{
N\exp\left[\frac{1}{N}\sum_{i=1}^{N} \log \left(1+\frac{\beta_{k,1}}{N} f_{k}(X_{i})\right)  - \frac{\log\frac{2m}{\varepsilon}}{N} \right]
-N}{D_{k}\beta_{k,1}}
 .$$
\end{thm}

Before we give the proof, let us see why this theorem really improves theorem \ref{firstbound}.
Let us choose put:
$$ V_{k} = P\left\{\left[f_{k}(X)-P\left(f_{k}(X)\right)\right]^{2} \right\}$$
and:
$$ \beta_{k,1} = \beta_{k,2} = \sqrt{\frac{N\log\frac{2m}{\varepsilon}}{V_{k}}} .$$
Then we obtain:
$$ \alpha^{\inf}_{k}(\varepsilon,\beta_{k,1}) = \hat{\alpha}_{k}-\sqrt{\frac{2V_{k}\log\frac{2m}{\varepsilon}}{N}} +
\mathcal{O}_{P} \left(\frac{\log\frac{2m}{\varepsilon}}{N} \right) $$
and:
$$ \alpha^{\sup}_{k}(\varepsilon,\beta_{k,2}) = \hat{\alpha}_{k}+\sqrt{\frac{2V_{k}\log\frac{2m}{\varepsilon}}{N}} +
\mathcal{O}_{P}\left(\frac{\log\frac{2m}{\varepsilon}}{N}\right).$$
So, the first order term for $d^{2}(\hat{\alpha}_{k}f_{k},\overline{\alpha}_{k}f_{k})$ is:
$$ \frac{2V_{k}\log\frac{2m}{\varepsilon}}{N} ,$$
there is an improvement by a factor $4$ when we compare this bound to theorem \ref{firstbound}.

Remark that this particular choice for $\beta_{k,1}$ and $\beta_{k,2}$ is valid as soon as:
$$ \sqrt{\frac{N\log\frac{2m}{\varepsilon}}{V_{k}}} < \frac{N}{\sqrt{C_{k}D_{k}}} $$
or equivalently as soon as $N$ is greater than
$$\frac{C_{k}D_{k}\log\frac{2m}{\varepsilon}}{V_{k}}.$$

In practice, however, this particular $\beta_{k,1}$ and $\beta_{k,2}$ are unknown.
We can use the following procedure (see Catoni \cite{Classif}).
We choose a value $a>1$ and:
$$ B = \left\{a^{l} ,0\leq l \leq \left\lfloor\frac{\log \frac{N}{\sqrt{C_{k}D_{k}}}}{\log a}\right\rfloor -1\right\}.$$
By taking a union bound over all possibles values of $B$, with:
$$|B| \leq \frac{\log \frac{N}{\sqrt{C_{k}D_{k}}}}{\log a}$$
we obtain the following corollary.

\begin{cor}
Under condition $\mathcal{H}(+\infty)$, for any $a>1$,
for any $\varepsilon>0$, with $P^{\otimes N}$-probability at least $1-\varepsilon$
we have:
$$
 \sup_{\beta\in B} \alpha^{\inf}_{k}\left(\frac{\varepsilon\log a}{\log N-\frac{1}{2}\log C_{k}D_{k} },\beta\right)
\leq
\overline{\alpha}_{k}
\leq
\inf_{\beta\in B} \alpha^{\sup}_{k}\left(\frac{\varepsilon\log a}{\log N -\frac{1}{2}\log C_{k}D_{k}},\beta\right)
,
$$
with:
$$ B = \left\{a^{l} ,0\leq l \leq \left\lfloor\frac{\log \frac{N}{\sqrt{C_{k}D_{k}}}}{\log a}\right\rfloor -1\right\}.$$
\end{cor}

Note that the price to pay for the optimization with respect to $\beta_{k,1}$ and $\beta_{k,2}$
was just a $\log\log N$ factor.

\begin{proof}[Proof of the theorem.]
The technique used in the proof is due to Catoni \cite{manuscrit}.
Let us choose $k\in\{1,...,m\}$, and:
$$\beta\in \left(0,\frac{N}{\sqrt{C_{k}D_{k}}}\right).$$
We have, for any $\eta\in\mathds{R}$:
$$
P^{\otimes N} \exp\left\{\sum_{i=1}^{N}\log \left(1-\frac{\beta}{N} f_{k}(X_{i})\right) -\eta \right\}
\leq
\exp\Biggl\{N \log\Bigl(1-\frac{\beta}{N} P\left[f_{k}(X)\right]\Bigr)-\eta\Biggr\}.
$$
Let us choose:
$$ \eta = \log\frac{2m}{\varepsilon} + N \log\Bigl(1-\frac{\beta}{N} P\left[f_{k}(X)\right]\Bigr) .$$
We obtain:
$$
P^{\otimes N} \exp\left\{\sum_{i=1}^{N}\log \left(1-\frac{\beta}{N} f_{k}(X_{i})\right) -
\log\frac{2m}{\varepsilon} - N \log\Bigl(1-\frac{\beta}{N} P\left[f_{k}(X)\right]\Bigr)
 \right\}
\leq
\frac{\varepsilon}{2m}
,
$$
and so:
$$
P^{\otimes N} \left\{\sum_{i=1}^{N}\log \left(1-\frac{\beta}{N} f_{k}(X_{i})\right) \geq
\log\frac{2m}{\varepsilon} + N \log\Bigl(1-\frac{\beta}{N} P\left[f_{k}(X)\right]\Bigr)
 \right\}
\leq \frac{\varepsilon}{2m},
$$
that becomes:
$$
P^{\otimes N} \left\{
P\left[f_{k}(X)\right]
\geq
\frac{N}{\beta}\left[
1-\exp\left(\frac{1}{N}\sum_{i=1}^{N} \log \left(1-\frac{\beta}{N} f_{k}(X_{i})\right) - \frac{\log\frac{2m}{\varepsilon}}{N} \right]
\right)
 \right\}
\leq \frac{\varepsilon}{2m}.
$$
We apply the same technique to:
$$
P^{\otimes N} \exp\left\{\sum_{i=1}^{N}\log \left(1+\frac{\beta'}{N} f_{k}(X_{i})\right) -\eta \right\}
\leq
\exp\Biggl\{N \log\Bigl(1+\frac{\beta'}{N} P\left[f_{k}(X)\right]\Bigr)-\eta\Biggr\}
$$
to obtain the upper bound.
We combine both result by a union bound argument.
\end{proof}

\subsection{A generalization to data-dependent basis functions}

We now extend the previous method to the case where the family  $(f_{1},...,f_{m})$ is allowed to be data-dependant, in a particular sense.
This subsection requires some modifications of the notations of section \ref{em}.

\begin{dfn}
For any $m'\in\mathds{N}^{*}$ we define a function $\Theta_{m'} : \mathcal{X} \rightarrow \left(\mathcal{L}^{2}\right)^{m'}$.
For any $i\in\{1,...,N\}$ we put:
$$ \Theta_{m'}(X_{i}) = \left(f_{i,1},...,f_{i,m'}\right) .$$
Finally, consider the family of functions:
$$ \left(f_{1},...,f_{m}\right) = \left(f_{1,1},...,f_{1,m'},...,f_{N,1},...,f_{N,m'}\right). $$
So we have $m=m'N$ (of course, $m'$ is allowed to depend on $N$).
Let us take, for any $i\in\{1,...,N\}$:
$$ P_{i} (.) = P^{\otimes N}(.|X_{i}) .$$
We put, for any $(i,k)\in\{1,...,N\}\times\{1,...,m'\}$:
$$ D_{i,k} = \int_{\mathcal{X}} f_{i,k}(x)^{2} \lambda(dx) ,$$
and we still assume that condition $\mathcal{H}(\infty)$ is satisfied, that means here that we have
known constants $ C_{i,k} = c_{i,k} $ such that:
$$ \forall x\in\mathcal{X}, \quad \left|f_{i,k}(x)\right| \leq \sqrt{C_{i,k}D_{i,k}} .$$
Finally, we put:
$$ \overline{\alpha}_{i,k} = \arg\min_{\alpha\in\mathds{R}} d^{2}(\alpha f_{i,k},f) .$$
\end{dfn}

Let us choose $(i,k)\in\{1,...,N\}\times\{1,...,m'\}$.
Using Seeger's idea, we follow the preceding proof, replacing $P^{\otimes N}$ by $P_{i}$, and using
the $N-1$ random variables:
$$ \Bigl(f_{i,k}(X_{j})\Bigr)_{
\tiny{
\begin{array}{l}
j\in\{1,...,N\}
\\
j\neq i
\end{array}
}
} $$
with
$$\eta=\log\frac{2m'N}{\varepsilon}+(N-1)\log\left(1-\frac{\beta}{N-1}P\left[f_{i,k}(X)\right]\right) $$
and we obtain:
\begin{multline*}
P_{i} \exp\Biggl\{\sum_{j\neq i}\log \left(1-\frac{\beta}{N-1} f_{i,k}(X_{j})\right) -
\log\frac{2m'N}{\varepsilon}
\\
- (N-1) \log\Bigl(1-\frac{\beta}{N-1} P\left[f_{i,k}(X)\right]\Bigr)
 \Biggr\}
\leq
\frac{\varepsilon}{2m'N}.
\end{multline*}
Note that for any random variable $H$ that is a function of the $X_{i}$:
$$ P^{\otimes N} P_{i} H = P^{\otimes N} H .$$
So we conclude exactly in the same way than for the previous theorem and we obtain the following result.

\begin{thm}
For any $\varepsilon>0$,
for any $\beta_{i,k,1},\beta_{i,k,2}$ such that:
$$  0 < \beta_{i,k,j} < \frac{N-1}{\sqrt{C_{i,k}D_{i,k}}} , \quad j\in\{1,2\} ,$$
with $P^{\otimes N}$-probability at least $1-\varepsilon$,
for any $i\in\{1,...,N\}$ and $k\in\{1,...,m\}$
we have:
$$ \tilde{\alpha}^{\inf}_{k}(\varepsilon,\beta_{i,k,1}) \leq \overline{\alpha}_{k} \leq \tilde{\alpha}^{\sup}_{k}(\varepsilon,\beta_{i,k,2}) $$
with:
\begin{multline*}
 \tilde{\alpha}^{\sup}_{k}(\varepsilon,\beta_{i,k,2})
\\  =
\frac{
N-1-(N-1)\exp\left[\frac{1}{N-1}\sum_{j\neq i} \log \left(1-\frac{\beta_{i,k,2}}{N-1} f_{i,k}(X_{j})\right)-\frac{\log\frac{2m'N}{\varepsilon}}{N-1}\right]
}{D_{i,k}\beta_{i,k,2}}
\end{multline*}
and:
\begin{multline*}
\tilde{\alpha}^{\inf}_{k}(\varepsilon,\beta_{i,k,1}) 
\\  =
\frac{
(N-1)\exp\left[\frac{1}{N-1}\sum_{j\neq i} \log \left(1+\frac{\beta_{i,k,1}}{N-1} f_{i,k}(X_{j})\right)- \frac{\log\frac{2m'N}{\varepsilon}}{N-1} \right]
-N+1}{D_{i,k}\beta_{i,k,1}}
 .
\end{multline*}
\end{thm}

\begin{exm}[Support Vector Machines]
\label{exSVM}
Actually, SVM were firstly introduced by Guyon, Boser and Vapnik \cite{SVM_FIRST} in the context of classification,
but the method was extended by Vapnik \cite{Vapnik} to the context of least square regression estimation and of
density estimation. The idea is to generalize the kernel estimator to the case where $\mathcal{X}$ is of large dimension,
and so we cannot use a grid like we did in the $[0,1]$ case.
Let us choose a function:
\begin{align*}
K: \mathcal{X}^{2} & \rightarrow\mathds{R} \\
              (x,x') & \mapsto K(x,x').
\end{align*}
We take $m'=1$ and:
$$ \Theta_{1}(x) = \left(K(x,.)\right)$$
then the obtained estimator has the form of a SVM:
$$ \hat{f}(x) = \sum_{i=1}^{N} \tilde{\alpha}_{i}K(X_{i},x) $$
where the set of $i$ such that $\tilde{\alpha}_{i}\neq 0$ is expected to be small.
Note that we do not need to assume that $K(.,.)$ is a Mercer's kernel as usual with SVM.
Moreover, we can extend the method to the case where we have several kernels $K_{1},...,K_{m'}$ by taking:
$$ \Theta_{m'}(x) = \left(K_{1}(x,.),...,K_{m'}(x,.)\right) .$$
The estimator becomes:
$$ \hat{f}(x) = \sum_{j=1}^{m'} \sum_{i=1}^{N} \tilde{\alpha}_{i,j}K_{j}(X_{i},x) .$$
Note that a widely used kernel is the gaussian kernel; let $\delta(.,.)$ be a distance on $\mathcal{X}$ and
$\gamma_{1},...,\gamma_{m'}>0$
then we put:
$$ K_{k}(x,x') = \exp\left(-\gamma_{k} \delta^{2}(x,x') \right). $$
For example, if $\mathcal{X}=\mathds{R}$ and $\lambda$ is the Lebesgue measure then
hypothesis $\mathcal{H}(\infty)$ is obviously satisfied with the gaussian kernel with
$$ C_{i,k} = c_{i,k} = \sqrt{\frac{\gamma_{k}}{\pi}} = \frac{1}{D_{i,k}} .$$
\end{exm}

\subsection{Back to the histogram}

In the case of the histogram, $f_{k}(.) = \mathds{1}_{A_{k}}(.)$ can take only two values: $0$ and $1$.
Remember that $D_{k}=\lambda(A_{k})$. So:
$$ \alpha^{\inf}_{k}(\varepsilon,\beta_{k,1}) =
\frac{N}{\lambda(A_{k}) \beta_{k,1}}
\left\{\left[\left(1+\frac{\beta_{k,1}}{N} \right)^{|\{i:X_{i}\in A_{k}\}|} \frac{\varepsilon}{2m}\right]^{\frac{1}{N}}-1\right\}
.
$$
Remember that, for any $x\geq 0$:
$$ (1+x)^{\gamma} \geq 1 + \gamma x + \frac{\gamma(\gamma-1)}{2} x^{2}  $$
and so:
$$
\alpha^{\inf}_{k}(\varepsilon,\beta_{k,1})
\geq
\hat{\alpha}_{k} \left(\frac{\varepsilon}{2m}\right)^{\frac{1}{N}}
\left[1-\frac{\beta_{k,1}(1-\hat{\alpha}_{k}D_{k})}{2N} \right]
- \frac{N}{D_{k}\beta_{k,1}} \left[1-\left(\frac{\varepsilon}{2m}\right)^{\frac{1}{N}}\right].
$$
Now, we take the grid:
$$ B = \left\{2^{l} ,0\leq l \leq \left\lfloor\frac{\log \frac{N}{\sqrt{D_{k}}}}{\log 2}\right\rfloor -1\right\}. $$
Remark that, for any $\beta$ in:
$$\left[1,\frac{N}{2\sqrt{D_{k}}}\right]$$
there is some $b\in B$ such that $\beta\leq b\leq 2\beta$, and so:
$$
\alpha^{\inf}_{k}(\varepsilon,b)
\geq
\hat{\alpha}_{k} \left(\frac{\varepsilon}{2m}\right)^{\frac{1}{N}}
\left[1-\frac{\beta_{k,1}(1-\hat{\alpha}_{k}D_{k})}{2N} \right]
- \frac{N}{D_{k}2\beta_{k,1}} \left[1-\left(\frac{\varepsilon}{2m}\right)^{\frac{1}{N}}\right].
$$
This allows us to choose whatever value for $\beta_{k,1}$ in
$$\left[1,\frac{N}{2\sqrt{D_{k}}}\right].$$
Let us choose:
$$ \beta_{k,1} = \sqrt{\frac{N^{2} \left[\left(\frac{\varepsilon}{2m}\right)^{\frac{-1}{N}}-1\right]}{\hat{\alpha}_{k}D_{k}(1-\hat{\alpha}_{k}D_{k})}} $$
that is allowed for $N$ large enough.
So we have:
$$
\alpha^{\inf}_{k}(\varepsilon,\beta_{k,1})
\geq
\hat{\alpha}_{k} \left(\frac{\varepsilon}{2m}\right)^{\frac{1}{N}}
- \sqrt{\hat{\alpha}_{k}D_{k}(1-\hat{\alpha}_{k}D_{k}) \left[\left(\frac{\varepsilon}{2m}\right)^{\frac{-1}{N}}-1\right] }.
$$
With the union bound term (over the grid $B$) we obtain:
\begin{multline*}
\alpha^{\inf}_{k}\left(\frac{\varepsilon \log 2}{\log \frac{N}{\sqrt{D_{k}}}},\beta_{k,1}\right)
\\
\geq
\hat{\alpha}_{k} \left(\frac{\varepsilon\log 2}{2m\log \frac{N}{\sqrt{D_{k}}}}\right)^{\frac{1}{N}}
- \sqrt{\hat{\alpha}_{k}D_{k}(1-\hat{\alpha}_{k}D_{k}) \left[\left(\frac{\varepsilon\log 2}{2m\log\frac{N}{\sqrt{D_{k}}}}\right)^{\frac{-1}{N}}-1\right]}
\\
=
\hat{\alpha}_{k}
- \sqrt{\frac{\hat{\alpha}_{k}D_{k}(1-\hat{\alpha}_{k}D_{k}) \log\frac{2m\log \frac{N}{\sqrt{D_{k}}}}{\varepsilon \log 2} }{N}}
+ \mathcal{O}\left(\frac{\log\frac{m\log N}{\varepsilon}}{N}\right),
\end{multline*}
remark that we have this time the "real" variance term of $\mathds{1}_{A_{k}}(X)$:
$$ \hat{\alpha}_{k}D_{k}(1-\hat{\alpha}_{k}D_{k}) = \frac{|\{i:X_{i}\in A_{k}\}|}{N}\left(1-\frac{|\{i:X_{i}\in A_{k}\}|}{N}\right) .$$

\subsection{Another simple example: the Haar basis}

\label{Haar}

Let us assume that $\mathcal{X} = [0,1]$.
Let $(\varphi,\psi)$ be a father wavelet and the associated mother wavelet, and:
$$\psi_{j,k}(x)=\psi(2^{j}x+k)$$
for $k\in\{0,...,2^{j}-1\}=S_{j}$ (note that the wavelet basis is non-normalized here).
Here, we use the Haar wavelets, with:
\begin{align*}
\varphi(x) & = \mathds{1}_{[0,1]}(x) \\
\psi(x) & = \mathds{1}_{\left[0,\frac{1}{2}\right]}(x) - \mathds{1}_{\left[\frac{1}{2},1\right]}(x)
.
\end{align*}
For the sake of simplicity, let us write:
$$\psi_{-1,k}(x)=\varphi(x)$$
for $k\in\{0\}=S_{-1}$.
By an obvious adaptation of our notations, let us put $\overline{\alpha}_{j,k}$ the coefficient associated to $\psi_{j,k}$:
$$ \overline{\alpha}_{j,k} = \frac{P\psi_{j,k}(X)}{\int \psi_{j,k}^{2} } = P\psi_{j,k}(X) ,$$
remark that condition $\mathcal{H}(\infty)$ is satisfied with
$D_{j,k}=2^{-j}$ and $C_{j,k}=1$.
In this particular setting, note that $ \overline{\alpha}_{-1,0} =1$ is known, so the associated confidence interval
is just $\{1\}$.
Moreover, here $\psi_{j,k}(X)$ can take only three values: $-1$, $0$ and $1$.
Let us put:
$$\overline{P} = \frac{1}{N} \sum_{i=1}^{N}\delta_{X_{i}}. $$
Remark that in this case we have:
\begin{multline*}
\frac{1}{N}\sum_{i=1}^{N} \log \left(1-\frac{\beta}{N} \psi_{j,k}(X_{i})\right)
\\
= \overline{P}(\psi_{j,k}(X)=1) \log \left(1-\frac{\beta}{N}\right)
+ \overline{P}(\psi_{j,k}(X)=-1) \log \left(1+\frac{\beta}{N}\right)
\\
= \frac{1}{2}\overline{P}\left[\psi_{j,k}(X)^{2}\right]\log\left(1-\frac{\beta^{2}}{N^{2}}\right)
+ \frac{1}{2}\overline{P}\left[\psi_{j,k}(X)\right]\log\left(\frac{1-\frac{\beta}{N}}{1+\frac{\beta}{N}}\right).
\end{multline*}
So we have:
\begin{multline*}
\alpha^{\sup}_{j,k}(\varepsilon,\beta)
\\
=
\frac{
N-N\exp\left[
\frac{1}{2}\overline{P}\left[\psi_{j,k}(X)^{2}\right]\log\left(1-\frac{\beta^{2}}{N^{2}}\right)
- \frac{1}{2}\overline{P}\left[\psi_{j,k}(X)\right]\log\left(\frac{1+\frac{\beta}{N}}{1-\frac{\beta}{N}}\right)
 - \frac{\log\frac{2m}{\varepsilon}}{N} \right]
}{D_{k}\beta_{k,2}}
\end{multline*}
and:
\begin{multline*}
\alpha^{\inf}_{j,k}(\varepsilon,\beta)
\\ =
\frac{
N\exp\left[
\frac{1}{2}\overline{P}\left[\psi_{j,k}(X)^{2}\right]\log\left(1-\frac{\beta^{2}}{N^{2}}\right)
+ \frac{1}{2}\overline{P}\left[\psi_{j,k}(X)\right]\log\left(\frac{1+\frac{\beta}{N}}{1-\frac{\beta}{N}}\right)
- \frac{\log\frac{2m}{\varepsilon}}{N} \right]
-N}{D_{k}\beta_{k,1}}
 .
\end{multline*}

\section{Simulations}

\label{si}

\subsection{Description of the example}

We assume that we observe $X_{i}$
for $i\in\{1,...,N\}$ with $N=2^{10}=1024$,
where the variables $X_{i}\in[0,1]\subset\mathds{R}$ are i.i.d. from a distribution
with an unknown density $f$ with respect to the Lebesgue measure.
The goal is to estimate $f$.

Here, we will use three methods. The first estimation method will be a multiple kernel estimator obtained
by the algorithm described previously, the second one a thresholded wavelets estimate
also obtained by this algorithm, and we will compare both estimators to
a thresholded wavelet estimate as given by Donoho, Johnstone, Kerkyacharian and Picard \cite{Donoho}.

\subsection{The estimators}

\subsubsection{Hard-thresholded wavelet estimator}

We first use a classical hard-thresholded wavelet estimator.

In the case of the Haar basis (see subsection \ref{Haar}), we take:
$$\hat{\alpha}_{j,k} = 2^{j} \frac{1}{N}\sum_{j=1}^{N}\psi_{j,k}(X_{i}) .$$
For a given $\kappa\geq 0$ and $J\in\mathds{N}$, we take:
$$\tilde{f}_{J}(.)=\sum_{j=-1}^{J}\sum_{k\in S_{j}}\hat{\alpha}_{j,k}\mathds{1}(|\hat{\alpha}_{j,k}|\geq \kappa t_{j,N})
\psi_{j,k}(.) $$
where:
$$t_{j,N}=\sqrt{\frac{j}{N}} .$$

Actually, we must choose $J$ in such a way that:
$$2^{J}\sim t_{N}^{-1} .$$
Here, we choose $\kappa=0.7$ and $J=7$.

\subsubsection{Wavelet estimators with our algorithm}

We also use the same family of functions, and we apply our thresholding method, with bounds given in subsection \ref{Haar}.
So we take:
$$m =2^{J}= 128 .$$

We use an asymptotic version of our confidence intervals inspired by our theoretical confidence intervals:
$$ \overline{\alpha}_{j,k} \in \left[\hat{\alpha}_{j,k} \pm \sqrt{\frac{\log\frac{2m}{\varepsilon} V_{j,k}}{N}} \right] $$
where $V_{j,k}$ is the estimated variance of $\psi_{j,k}(X)$:
$$ V_{j,k} = \frac{1}{N}\sum_{i=1}^{N} \left[\psi_{j,k}(X_{i}) - \frac{1}{N}\sum_{h=1}^{N} \psi_{j,k}(X_{h})  \right]^{2} . $$

Let us remark that the union bound are always "pessimistic", and that
we use a union bound argument over all the $m$ models despite only a few of them
are effectively used in the estimator. So, we propose to actually use the individual confidence interval for each model, replacing:
the $ \log\frac{2m}{\varepsilon}$ by $\log\frac{2}{\varepsilon}$.

\subsubsection{Multliple estimator}

Finally, we use the kernel estimator described in section \ref{rc}, with function $K$:
$$ K_{j}(u,v) = \exp\left[-2^{2j}(u-v)^{2} \right] $$
with $n=N$ and $j\in\{1,...,h=6\}$.
We add the constant function $1$ to the family.

Here again we use the individuals confidence intervals, and the asymptotic version of this intervals.

\subsection{Experiments and results}

The simulations were realized with the R software \cite{R}.

For the experiments, we use the following functions $f$ that are some variations of the functions
used by Donoho and Johnstone for experiments on wavelets, for example in \cite{Donoho} (actually, these
functions were used as regression functions, so the modification was to add them a constant in order to ensure they take
nonnegative values):
\begin{align*}
Doppler(t) & = 1+2\sqrt{t(1-t)}\sin\frac{2\pi(1+v)}{t+v} \quad \text{ where } v=0.05 \\
HeaviSine(t) & = 1.5+ \frac{1}{4}\Bigl[4\sin 4\pi t - sgn(t-0.3)-sgn(0.72-t)\Bigr] \\
Blocks(t) & = 1.05+ \frac{1}{4}\sum_{i=1}^{11} c_{i} \mathds{1}_{(t_{i},+\infty)}(t)
\end{align*}
where $sgn(t)$ is the sign of $t$ (say $-1$ if $t\leq 0$ and $+1$ otherwise).
The values of the $c_{i}$ and $t_{i}$ are given in figure 1.

\begin{figure}
\caption{Values of $t_{i}$ and $c_{i}$ in the fonction $Blocks(.)$.}
\small{
\begin{tabular}{|p{0.5cm}|p{0.5cm}|p{0.5cm}|p{0.5cm}|p{0.5cm}|p{0.5cm}|p{0.5cm}|p{0.5cm}|p{0.5cm}|p{0.5cm}|p{0.5cm}|p{0.5cm}|}
\hline
$i$ & $1$ & $2$ & $3$ & $4$ & $5$ & $6$ & $7$ & $8$ & $9$ & $10$ & $11$ \\
\hline
$c_{i}$ & $4$ & $-5$ & $3$ & $-4$ & $5$ & $4.2$ & $-2.1$ & $4.3$ & $-3.1$ & $2.1$ & $-4.2$ \\
\hline
$t_{i}$ & $0.10$ & $0.13$ & $0.15$ & $0.23$ & $0.25$ & $0.40$ & $0.44$ & $0.65$ & $0.76$ & $0.78$ & $0.81$ \\
\hline
\end{tabular}
}
\end{figure}

We consider 3 experiments (for the three density functions), we choose $\varepsilon$=10\%, repeat each experiment 20 times;
the results are reported in figure 2. We also give some illustrations (figure 3, 4 and 5).

\begin{figure}
\caption{Results of the experiments. For each experiment, we give the mean distance of the estimator the the density ($d^{2}(.,f)$).}
\begin{tabular}{|p{1.7cm}|p{3cm}|p{3cm}|p{3cm}|}
\hline
Function $f(.)$ & standard thresholded wavelets & thresh. wav. with our method & multiple kernel \\
\hline $Doppler$   & 0.104 & 0.127 & 0.083 \\
\hline $HeaviSine$ & 0.071 & 0.066 & 0.040 \\
\hline $Blocks$    & 0.110 & 0.142 & 0.121 \\
\hline
\end{tabular}
\end{figure}

\begin{figure}[h]
\includegraphics[width=8cm, height=8cm] {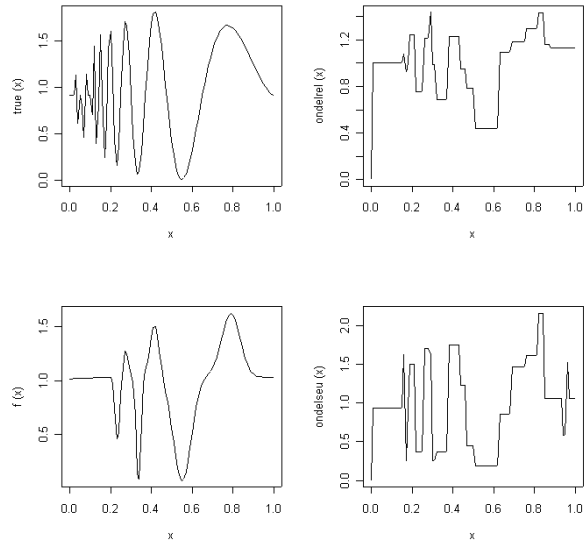}
\caption{Experiment 1, $f=Doppler$.
\small{Up-left: true regression function ($true$). Down-left: SVM ($f$).
Up-right: wavelet estimate with our algorithm ($ondelrel$).
Down-right: "classical" wavelet estimate ($ondelseu$).}}
\end{figure}

\begin{figure}[h]
\includegraphics[width=8cm, height=8cm] {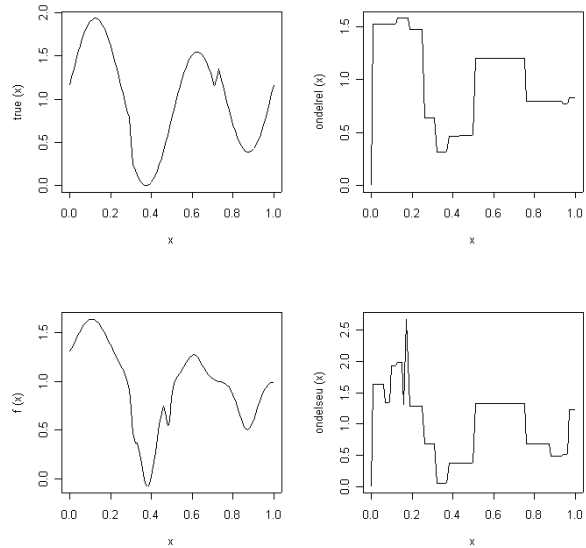}
\caption{Experiment 2, $f=HeaviSine$ and $\sigma=0.3$.}
\end{figure}

\begin{figure}[h]
\includegraphics[width=8cm, height=8cm] {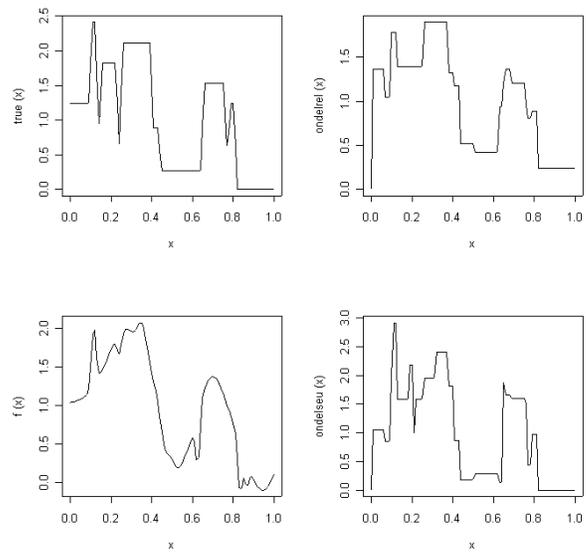}
\caption{Experiment 3, $f=Blocks$ and $\sigma=0.3$.}
\end{figure}

\end{document}